\documentclass[11pt]{article}%
\usepackage{amssymb}
\usepackage{times,nicefrac}
\usepackage{latexsym}
\usepackage{amsmath}
\usepackage{graphicx}
\usepackage{epsfig}
\usepackage{amsfonts}%
\providecommand{\U}[1]{\protect\rule{.1in}{.1in}}
\providecommand{\U}[1]{\protect\rule{.1in}{.1in}} \textwidth6.35in
\newcommand{\nc}{\newcommand}
\nc{\nn}{\nonumber} \nc{\eps}{\varepsilon} \nc{\la}{\lambda}
\nc{\wht}{\widehat} \nc{\ov}{\overline} \nc{\ds}{\displaystyle}
\nc{\ts}{\textstyle} \nc{\kro}{\left(} \nc{\kvo}{\left[}
\nc{\fio}{\left\{} \nc{\krz}{\right)} \nc{\kvz}{\right]}
\nc{\fiz}{\right\}} \nc{\noi}{\noindent}

\newtheorem{theorem}{Theorem}
\newtheorem{lemma}[theorem]{Lemma}

\newtheorem{corollary}[theorem]{Corollary}

\begin{document}

\title{Covering spheres with spheres}
\date{}
\author{Ilya Dumer\thanks{The author is with the University of California, Riverside,
USA (e-mail: dumer@ee.ucr.edu) }}
\maketitle

\begin{abstract}
Given a sphere of any radius $r$ in an $n$-dimensional Euclidean space, we
study the coverings of this sphere with\ solid spheres of radius one. Our goal
is to design a covering of the lowest \textit{covering density}, which defines
the average number of solid spheres covering a point in a bigger sphere. For
growing dimension $n,$ we design a covering that has covering density of order
$(n\ln n)/2$ for the full Euclidean space or for a sphere of any radius $r>1.$
This new upper bound reduces two times the asymptotic order of $n\ln n$
established in the classical Rogers bound.

\end{abstract}

\section{Introduction}

\textit{Spherical coverings. \ \ }Let $B_{r}^{n}(\mathbf{x})$ be
a\textit{\ ball} (solid sphere) of radius $r$ centered at some point
$\mathbf{x}=(x_{1},\ldots,x_{n})$ of an $n$-dimensional Euclidean space
$\mathbb{R}^{n}:$%
\[
B_{r}^{n}(\mathbf{x})\overset{\text{def}}{=}\left\{  \mathbf{z}\in
\mathbb{R}^{n}\mid\sum_{i=1}^{n}(z_{i}-x_{i})^{2}\leq r^{2}\right\}  .
\]
We also use a simpler notation $B_{r}^{n}$ if a ball is centered at the origin
$\mathbf{x}=0.$ For any subset $A\subseteq\mathbb{R}^{n}$, we say that a
subset \textsf{Cov}$(A,\varepsilon)\subseteq\mathbb{R}^{n}$ forms an
$\varepsilon$-\textit{covering} (an $\varepsilon$-net) of $A$ if $A$ is
contained in the union of the balls of radius $\varepsilon$ centered at points
$\mathbf{x}\in$\textsf{Cov}$(A,\varepsilon).$ In this case, we use notation
\[
\text{\textsf{Cov}}(A,\varepsilon):\;A\subseteq\bigcup_{\mathbf{x}%
\in\text{\textsf{Cov}}(A,\varepsilon)}B_{\varepsilon}^{n}(\mathbf{x}).
\]
By changing the scale in $\mathbb{R}^{n},$ we can always consider the rescaled
set $A/\epsilon$ and the new covering \textsf{Cov}$(A/\epsilon,1)$ with unit
balls $B_{1}^{n}(\mathbf{x}).$ Without loss of generality, below we consider
these (unit) coverings$.$ One of the classical problems is to obtain tight
bounds on the covering size $\left\vert \text{\textsf{Cov}}(B_{r}%
^{n},1)\right\vert $ for any ball $B_{r}^{n}$ of radius $r$ and dimension $n.$

Another related covering problem arises for a \textit{sphere}
\[
S_{r}^{n}\overset{\text{def}}{=}\left\{  \mathbf{z}\in\mathbb{R}^{n+1}\mid
\sum_{i=1}^{n+1}z_{i}^{2}=r^{2}\right\}  .
\]
Then a unit ball $B_{1}^{n+1}(\mathbf{x})$ intersects this sphere with a
spherical cap
\[
C_{r}^{n}(\rho,\mathbf{y})=S_{r}^{n}\cap B_{1}^{n+1}(\mathbf{x}),
\]
which has some center $\mathbf{y}\in S_{r}^{n}$, half-chord $\rho\leq1,$ and
the corresponding half-angle $\alpha=\arcsin\rho/r$. \ The biggest possible
cap $C_{r}^{n}(1,\mathbf{y})$ is obtained if the center $\mathbf{x}$ of the
corresponding ball $B_{1}^{n+1}(\mathbf{x})$ is centered at the distance
\begin{equation}
||\mathbf{x}||=\sqrt{r^{2}-1} \label{cen-x}%
\end{equation}
from the origin. To obtain a minimal covering, we shall consider the biggest
caps $C_{r}^{n}(1,\mathbf{y})$ assuming that all the centers $\mathbf{x}$
satisfy (\ref{cen-x})$.$\

\textit{Covering density.} Given a set $A\subseteq\mathbb{R}^{n},$ let
$\left|  A\right|  $ denote $n$-dimensional volume (Lebesque measure) of $A.$
We then consider any unit covering \textsf{Cov}$(A,1)$ and define minimum
covering density
\[
\vartheta(A)=\underset{\text{\textsf{Cov}}(A,1)}{\min}\;\sum_{\mathbf{x}%
\in\text{\textsf{Cov}}(A,1)}\frac{\left|  B_{1}^{n}(\mathbf{x})\cap A\right|
}{\left|  A\right|  }.
\]
Minimal coverings have been long studied for the spheres $S_{r}^{n}$ and the
balls $B_{r}^{n}.$ The celebrated Coxeter-Few-Rogers lower\ bound \cite{cox}
shows that for a sufficiently large ball $B_{r}^{n},$
\begin{equation}
\vartheta(S_{r}^{n})\geq c_{_{0}}n. \label{few}%
\end{equation}
Here and below $c_{i}$ denote some universal constants. A similar result also
holds for any sphere $S_{r}^{n}$ of radius $r\geq n$ (see Example 6.3 in
\cite{5}).

Various upper bounds on the minimum covering density are obtained for
$B_{r}^{n}$ and $S_{r}^{n}$ by Rogers in the classic papers \cite{1} and
\cite{3}. In particular, it follows from these papers that for a sufficiently
large radius $r,$ a ball $B_{r}^{n}$ and a sphere $S_{r}^{n}$ can be covered
with density
\begin{equation}
\vartheta\leq\left(  1+\frac{\ln\ln n}{\ln n}+\frac{5}{\ln n}\right)  n\ln n.
\label{d-ball}%
\end{equation}
Despite recent improvements obtained in \cite{5} and \cite{7}, respectively,
for spheres $S_{r}^{n}$ and balls\ $B_{r}^{n}$ of a relatively small radius
$r$, \ the Rogers bound (\ref{d-ball}) is still the best asymptotic bound
known for sufficiently large spheres, balls, and complete spaces
$\mathbb{R}^{n}$ of growing dimension $n.$

For a sphere $S_{r}^{n}$ of any dimension $n\geq3$ and an arbitrary radius
$r>1$, the best universal upper bound known to date is obtained in \cite{5}
(see Corollary 1.2 and Remark 5.1):
\begin{equation}
\vartheta_{\ast}\leq\left(  1+\frac{2}{\ln n}\right)  \left(  1+\frac{\ln\ln
n}{\ln n}+\frac{\sqrt{e}}{n\ln n}\right)  n\ln n. \label{d-sphe}%
\end{equation}
Our main result is presented in Theorem \ref{th:new}, which reduces about two
times the present upper bounds (\ref{d-ball}) and (\ref{d-sphe}) for
$n\rightarrow\infty.$

\begin{theorem}
\label{th:new}Unit balls can cover a sphere $S_{r}^{n}$ of any radius $r>1$
and any dimension $n\geq3$ with density%
\begin{equation}
\vartheta(S_{r}^{n})\leq\left(  \frac{1}{2}+\frac{3\ln\ln n}{\ln n}+\frac
{3}{\ln n}\right)  n\ln n.\label{d-sph}%
\end{equation}

\end{theorem}

The following corollary to Theorem \ref{th:new} ( see also \cite{dum1} ) shows
that the Rogers bound can also be reduced about two times for the coverings of
complete Euclidean spaces $\mathbb{R}^{n}.$

\begin{corollary}
\label{th:ball}For $n\rightarrow\infty,$ unit balls can cover the entire
Euclidean space $\mathbb{R}^{n}$ with density%
\begin{equation}
\vartheta(\mathbb{R}^{n})\leq\left(  \frac{1}{2}+o(1)\right)  n\ln n.
\label{d-spa}%
\end{equation}

\end{corollary}

\section{Preliminaries: embedded coverings\label{s:2}}

In this section, we obtain an estimate on $\vartheta(S_{r}^{n})$ that is
similar to { (\ref{d-sphe})} but uses a different technique. This technique of
embedded coverings will be substantially extended in Section \ref{s:3} to
improve the former bounds {(\ref{d-ball}) and (\ref{d-sphe})}. \ We will also
use most of our calculations performed in this section.

Consider a sphere $S_{r}^{n}$ of some dimension $n\geq3$ and radius $r>1.$ We
use notation $C(\rho,\mathbf{y})$ for a cap $C_{r}^{n}(\rho,\mathbf{y})$
whenever parameters $n$ and $r$ are fixed; we also use a shorter notation
$C(\rho)$ when a specific center $\mathbf{y}$\ is of no importance. In this
case, \textsf{Cov}$(\rho)$ will denote any covering of $S_{r}^{n}$ with
spherical caps $C(\rho).$ By definition, \ a covering \textsf{Cov}$(\rho)$ has
covering density
\[
\vartheta_{\rho}=\Omega_{\rho}\left\vert \text{\textsf{Cov}}(\rho)\right\vert
\]
where $\Omega_{\rho}$ is the fraction of the surface of the sphere $S_{r}^{n}$
covered by a cap $C\left(  \rho\right)  ,$
\[
\Omega_{\rho}\text{ }=\frac{\left\vert C\left(  \rho\right)  \right\vert
}{\left\vert S_{r}^{n}\right\vert }%
\]
We begin with two preliminary lemmas, which will simplify our calculations.
Let $f_{1}(x)$ and $f_{2}(x)$ be two positive differentiable functions$.$ We
say that $f_{1}(x)$ moderates $f_{2}(x)$ for $x\geq a$ if for all $x\geq a,$%
\[
\frac{f_{1}^{\prime}(x)}{f_{1}(x)}\geq\frac{f_{2}^{\prime}(x)}{f_{2}(x)}.
\]

\begin{lemma}
\label{lm:mod}Consider $m$\ functions $f_{i}(x)$ such that $f_{1}(x)$
moderates each function $f_{i}(x),$ $i\geq2,$\ for $x\geq a.$ Then inequality
$f_{1}(x)\geq\sum_{i=2}^{m}f_{i}(x)$ holds for any $x\geq a$ if it is valid
for $x=a.$
\end{lemma}

\noindent\textit{Proof. }Note that$\ f_{i}(x)=f_{i}(a)\exp\left\{
s_{i}(x)\right\}  ,$ where $s_{i}(x)\overset{\text{def}}{=}\int_{a}^{x}%
\frac{f_{i}^{\prime}(t)}{f_{i}(t)}dt.$ Also, $s_{i}(x)\leq s_{1}(x)$ for all
$i\geq2.$ Therefore,
\[
f_{1}(x)\geq\sum_{i=2}^{m}f_{i}(a)\exp\left\{  s_{1}(x)\right\}  \geq
\sum_{i=2}^{m}f_{i}(a)\exp\left\{  s_{i}(x)\right\}  =\sum_{i=2}^{m}f_{i}(x),
\]
which completes the proof.{\hfill} {\ }${\square}$\smallskip\smallskip

Let $\rho=1-\varepsilon.$ We first estimate the volumes of the caps
$C(\varepsilon,\mathbf{x})$ and $C(\rho,\mathbf{x})$ in relation to the volume
$\Omega_{1}=1/2$ of a hemisphere $C(1).$ Let $k_{n}$ be the volume of the unit
Euclidean $n$-ball $B_{n}.$

\begin{lemma}
\label{lm:epsi}The caps $C(\varepsilon,\mathbf{x})$ and $C(\rho,\mathbf{x})$
in a sphere $S_{1}^{n}$ have volumes%
\begin{align}
\Omega_{\varepsilon}  &  \geq\varepsilon^{n}/\left(  3\sqrt{1-\varepsilon^{2}%
}\sqrt{(n+1)2\pi}\right) \label{p-a}\\
\Omega_{\rho}  &  \geq\frac{1}{2}-\sqrt{\frac{n\varepsilon}{\pi}%
},\;\varepsilon<\pi/4n \label{expan2}%
\end{align}

\end{lemma}

\noindent\textit{Proof.} Inequality {(\ref{p-a}) follows from } \cite{5},
Lemma 3.1. Next, we prove (\ref{expan2}). Let $\alpha=\arcsin\rho.$ Then a cap
$C(\rho,\mathbf{x})$ have volume
\[
\Omega_{\rho}=nk_{n}\int_{0}^{\alpha}\sin^{n}\beta\,d\beta\geq\Omega
_{1}-nk_{n}\left(  \pi/2-\alpha\right)
\]
Note that $\pi/2-\alpha\leq\sqrt{2\varepsilon}.$ Indeed,
\[
\sin(\pi/2-\alpha)=\cos\alpha=\sqrt{2\varepsilon-\varepsilon^{2}}\leq
\sqrt{2\varepsilon}(1-\varepsilon/4)
\]
On the other hand, for any $\varepsilon<3/4$
\[
\sin(\sqrt{2\varepsilon})\geq\sqrt{2\varepsilon}(1-\varepsilon^{2}/3)\geq
\sqrt{2\varepsilon}(1-\varepsilon/4)
\]
Here we used inequality $\sin x\geq x(1-x^{2}/6)$ for any $x\in(0,\pi/2).$
Thus, we obtain (\ref{expan2}) since%
\begin{align*}
\Omega_{1}  &  \geq k_{n}\sqrt{n\pi/2}\\
\Omega_{1}-\Omega_{\rho}  &  \leq nk_{n}\sqrt{2\varepsilon}\leq\Omega_{1}%
\sqrt{4n\varepsilon/\pi}%
\end{align*}
{\hfill}${\square}$\smallskip

\textit{An embedded algorithm. }We employ the following parameters
\begin{equation}
\varepsilon=\frac{1}{2n\ln^{2}n},\;\rho=1-\varepsilon,\;\lambda=1+\frac
{2\ln\ln n}{\ln n} \label{cap2}%
\end{equation}
For $n\geq20,$ these parameters {simplify bounds } {(\ref{p-a}) } as follows {
}%
\begin{equation}%
\begin{tabular}
[c]{l}%
$\Omega_{\varepsilon}\geq\varepsilon^{n}/\left(  8\sqrt{n}\right)  $\\
$\Omega_{\rho}\geq\theta_{n}=\frac{1}{2}-\frac{1}{\ln n}\sqrt{\frac{1}{2\pi}}$%
\end{tabular}
\ \label{p01}%
\end{equation}
Here the first bound for $\Omega_{\varepsilon}$ is verified numerically.

To design a covering \textsf{Cov}$(1)$ of the sphere $S_{r}^{n},$ we first
randomly choose $N$ points $\mathbf{y\in}S_{r}^{n},$ where%
\begin{equation}
\frac{\lambda n\ln n}{\Omega_{\rho}}-1<N\leq\frac{\lambda n\ln n}{\Omega
_{\rho}} \label{p3}%
\end{equation}
Given the set $\left\{  C(\rho,\mathbf{y})\right\}  $ of $N$ caps, we then
consider another covering
\[
\text{\textsf{Cov}}(\varepsilon):\;S_{r}^{n}\subseteq\bigcup_{\mathbf{u}%
\in\text{\textsf{Cov}}(\varepsilon)}C(\varepsilon,\mathbf{u})
\]
with smaller caps $C(\varepsilon,\mathbf{u})$. Then we take all centers
$\mathbf{u}^{\prime}\in\,$\textsf{Cov}$(\varepsilon)$ that are left uncovered
by the set $\left\{  C(\rho,\mathbf{y})\right\}  \ \ $and form the extended
set $\left\{  \mathbf{x}\right\}  =\{\mathbf{y\}\cup}\{\mathbf{u}^{\prime
}\mathbf{\}}.$ This set covers the entire set \textsf{Cov}$(\varepsilon).$ By
expanding the caps $C(\rho,\mathbf{x})$ to the caps $C(1,\mathbf{x),}$we
obtain a unit covering \
\[
\text{\textsf{Cov}}(1):\;S_{r}^{n}\subseteq\bigcup_{\mathbf{x}\in\left\{
\mathbf{x}\right\}  }C(1,\mathbf{x}).
\]
The following lemma yields slightly larger residual terms for density
$\vartheta(S_{r}^{n})$ than those obtained in {(\ref{d-sphe}); however, it
will allow us to further improve estimates in Section }\ref{s:3}.

\begin{lemma}
\label{lm:basic}For $n\geq20,$ covering $\left\{  \mathbf{x}\right\}  $ of a
sphere $S_{r}^{n}$ with unit caps $C(1,\mathbf{x)}$ has density
\begin{equation}
\vartheta(S_{r}^{n})\leq\left(  1+\frac{2\ln\ln n}{\ln n}+\frac{2}{\ln
n}\right)  n\ln n. \label{est0}%
\end{equation}

\end{lemma}

\noindent\textit{Proof. }Any point $\mathbf{u}$ is covered by some cap
$C(\rho,\mathbf{y})$ with probability $\Omega_{\rho}.$ Let $N^{\prime}$ be the
expected number of centers $\mathbf{u}^{\prime}$ in \textsf{Cov}%
$(\varepsilon)$ that are left uncovered after $N$ random trials. Then the
lower bound of (\ref{p3}) gives
\[
N^{\prime}=\left(  1-\Omega_{\rho}\right)  ^{N}\cdot\left\vert
\text{\textsf{Cov}}(\varepsilon)\right\vert \leq2\left(  1-\Omega_{\rho
}\right)  ^{\left(  \lambda n\ln n\right)  /\Omega_{\rho}}%
\]
Note that $(1-x)^{1/x}$ declines with $x\in(0,1).$ We then take the minimum
value $\Omega_{\rho}=\theta_{n}$ of (\ref{p01}). This gives the upper bound
\begin{align*}
\left(  1-\Omega_{\rho}\right)  ^{N}  &  \leq2\left(  1-\theta_{n}\right)
^{\left(  \lambda n\ln n\right)  /\theta_{n}}\leq2e^{nt_{n}},\\
t_{n}  &  =(\ln n+2\ln\ln n)\frac{\ln\left(  1-\theta_{n}\right)  }{\theta
_{n}}%
\end{align*}
Finally, inequality (\ref{p01}) gives
\begin{equation}
\left\vert \text{\textsf{Cov}}(\varepsilon)\right\vert =\vartheta
_{\varepsilon}/\Omega_{\varepsilon}\leq\left(  \vartheta_{\varepsilon}%
/\Omega_{1}\right)  \left[  2n\ln^{2}n\right]  ^{n}8\sqrt{n} \label{delta}%
\end{equation}
and%
\begin{align}
N^{\prime}  &  \leq2e^{nt_{n}}\left\vert \text{\textsf{Cov}}(\varepsilon
)\right\vert \leq2e^{nd_{n}}\vartheta_{\varepsilon}/\Omega_{1}, \label{delta1}%
\\
d_{n}  &  =t_{n}+\frac{\ln8}{n}+\frac{\ln n}{2n}+\ln2+\ln\left(  n\ln
^{2}n\right) \nonumber
\end{align}
Note that function $d_{n}$ has moderating term $t_{n}$. Then straightforward
calculations give $d_{n}\leq d_{20}<-0.4.$

Consider a covering $\left\{  \mathbf{x}\right\}  $ with caps $C(1,\mathbf{x)}%
$ that has size at most $N+N^{\prime}.$ Then (\ref{p3}) and (\ref{expan2})
give the covering density of $\left\{  \mathbf{x}\right\}  :$
\begin{align}
\vartheta_{1}  &  =\Omega_{1}(N+N^{\prime})\leq\left(  \lambda n\ln n\right)
\Omega_{1}/\Omega_{\rho}+\vartheta_{\varepsilon}e^{-n/3}\nonumber\\
&  \leq\left(  \lambda n\ln n\right)  /\left(  1-\sqrt{2/\pi}\ln^{-1}n\right)
+e^{-n/3}\vartheta_{\varepsilon} \label{eps5}%
\end{align}
For any given $n,$ bound (\ref{eps5}) only depends on the density
$\vartheta_{\varepsilon}$. Next, we can change the scale in $\mathbb{R}^{n+1}$
\ and replace a covering \textsf{Cov}$(1)$ of a sphere $S_{r/\varepsilon
}^{\,n}$ with the covering \textsf{Cov}$(\varepsilon)$ of the sphere
$S_{r}^{\,n}.$ This rescaling shows that we can replace $\vartheta
_{\varepsilon}$ in (\ref{eps5}) with any known density $\vartheta_{1}$. Thus,
this iteration process yields the upper bound
\begin{equation}
\vartheta_{1}\leq\frac{\lambda n\ln n}{\left(  1-\sqrt{2/\pi}\ln^{-1}n\right)
\left(  1-e^{-n/3}\right)  }. \label{eps7}%
\end{equation}
which we replace with a weaker bound (\ref{est0}). Here we again use \ Lemma
\ref{lm:mod} and verify that{\ the estimate }(\ref{est0}) exceeds{ the
estimate }(\ref{eps7}) {for }$n=20$ and moderates it for larger $n$.{\hfill}
{\ }${\square}$

\section{New covering algorithm for a sphere $S_{r}^{n}\label{s:3}$}

\textit{Covering design. }In this section, we obtain a covering of the sphere
$S_{r}^{n}$ with asymptotic density $\left(  n\ln n\right)  /2$. The new
design will use \ both the former covering \textsf{Cov}$(\varepsilon)$ and
another covering \textsf{Cov}$(\mu)$ with a larger radius $\mu$ that has
asymptotic order of $n^{-1/2}.$ \ Namely, we use parameters
\begin{align}
\varepsilon &  =\left(  2n\ln^{2}n\right)  ^{-1}\underset{}{,}\text{ \ }%
\rho=1-\varepsilon\nonumber\\
\beta &  =\underset{}{\frac{1}{2}}+\frac{3\ln\ln n}{\ln n},\text{ \ }%
\lambda=\beta+\frac{1}{2\ln n}\nonumber\\
\mu &  =n^{-\beta}/4=\underset{}{1}/\left(  4\sqrt{n}\ln^{3}n\right)
\nonumber\\
d &  =1-2\varepsilon-\mu^{2}=1-\frac{1}{n\ln^{2}n}-\frac{1}{16n\ln^{6}n}_{{}%
}\label{d0}%
\end{align}
\smallskip\smallskip\smallskip and proceed as follows.

A. Let a sphere $S_{r}^{n}$ be covered with two different coverings
\textsf{Cov}$(\mu)$ and \textsf{Cov}$(\varepsilon):$
\[
\text{\textsf{Cov}}(\mu):\;S_{r}^{n}\subseteq\bigcup_{\mathbf{z}%
\in\text{\textsf{Cov}}(\mu)}C(\mu,\mathbf{z}),
\]%
\[
\text{\textsf{Cov}}(\varepsilon):\;S_{r}^{n}\subseteq\bigcup_{\mathbf{u}%
\in\text{\textsf{Cov}}(\varepsilon)}C(\varepsilon,\mathbf{u}).
\]
We assume that both coverings have the former density $\vartheta_{\ast}$ of
(\ref{est0}) or less$.$

B. Randomly choose $N$ points $\mathbf{y}\in S_{r}^{n}$ and consider the
corresponding spherical caps $C(\rho,\mathbf{y}),$ where
\begin{equation}
N=\left\lfloor \left(  \lambda n\ln n\right)  /\Omega_{d}\right\rfloor
\label{n1}%
\end{equation}

C. Let $C(\mu,\mathbf{\bar{z}})$ be any cap in \textsf{Cov}$(\mu)$ that
contains at least one center $\mathbf{u}\in$\textsf{Cov}$(\varepsilon)$ \ not
covered by the $\rho$-caps$.$ We consider all such centers $\mathbf{\bar{z}}$
and form the joint set $\left\{  \mathbf{x}\right\}  =\left\{  \mathbf{y}%
\right\}  \cup\{\mathbf{\bar{z}}\}.$ This set covers \textsf{Cov}%
$(\varepsilon)$ with $\rho$-caps and therefore forms the required covering, by
extension to the caps $C(1,\mathbf{x}):$
\[
\text{\textsf{Cov}}(1):\;S_{r}^{n}\subseteq\bigcup_{\mathbf{x}\in\left\{
\mathbf{x}\right\}  }C(1,\mathbf{x}).
\]
We now proceed with preliminary discussion, which outlines the main steps of
the proof. \medskip

\textit{Outline of the proof. } Let us first assume that\ we keep the design
of Section \ref{s:2} but apply it to the new covering \textsf{Cov}$(\mu)$
instead of \textsf{Cov}$(\varepsilon).$ This will require taking $\rho=1-\mu$
to cover the centers of the caps $C(\mu,\mathbf{z})$\ and then expanding
$\rho$ to 1 to cover the whole $\mu$-caps. Contrary to our former choice of
$\rho=1-\varepsilon$ in (\ref{cap2}), \ it can be proven that this expansion
causes the covering density to grow exponentially in $n$. To circumvent this
problem, we keep $\rho=1-\varepsilon$ in (\ref{d0}) but change our design as
follows.\smallskip\smallskip

1. Given any cap $C(\mu,\mathbf{z}),$ we say that a cap $C(\rho,\mathbf{y})$
is $t$-close if $\mathbf{y}$ falls within distance $t<\rho$ to $\mathbf{z.}$
In our proof, we refine the selection of the caps $C(\rho,\mathbf{y})$ and
count only $d$-close caps$,$ instead of the $\rho$-close caps considered in
Section \ref{s:2}. \ It is easy to verify that distance $d$ of (\ref{d0}) is
so close to $\rho$ \ that
\begin{equation}
\Omega_{\rho}/\Omega_{d}\rightarrow1,\quad n\rightarrow\infty. \label{n2}%
\end{equation}
For this reason, counting only $d$-close caps instead of the former $\rho
$-close caps will carry no overhead to the covering size (\ref{n1}). \

2. On the other hand, we will show in Lemma {\ }\ref{lm:lower} \ that the
$\mu$-cap becomes almost completely covered by a cap $C(\rho,\mathbf{y})$ when
the latter becomes $d$-close instead of $\rho$-close$.$ Namely, only a
vanishing fraction $\omega<\exp\left(  -2\ln^{2}n\right)  $ of \ a $\mu$-cap
is left uncovered in this case.

3. We shall also use the fact that a typical $\mu$-cap is covered by multiple
$d$-close caps. According to (\ref{n1}), the average number $\Omega_{d}N$ of
these $d$-close caps has the exact order of $\lambda n\ln n:$
\begin{equation}
\lambda n\ln n-\Omega_{d}<\Omega_{d}N\leq\lambda n\ln n \label{lam1}%
\end{equation}
In our proof, we first define insufficiently covered $\mu$-caps. Namely, we
call a cap $C(\mu,\mathbf{z}^{\prime})$ non-saturated if it has only $s$ or
fewer $d$-close $\rho$-caps, where $s$ has a lower order,
\begin{equation}
s=\left\lfloor n/q\right\rfloor ,\quad q=6\ln\ln n \label{s}%
\end{equation}
This choice of $s$ will achieve two goals. \

4. We prove in Lemma \ref{lm:bad} that non-saturated caps typically form a
very small fraction among all $\mu$-caps. This fraction has the order below
$\exp\left[  -\lambda n\ln n\right]  $. On the other hand, it is easy to see
that the quantity%
\[
\left\vert \text{\textsf{Cov}}(\mu)\right\vert \leq{\vartheta}\left(
S_{r}^{n}\right)  /{\Omega_{\mu}}%
\]
exceeds $N$ by the exponential factor $\Omega_{d}/\Omega_{\mu}\sim\exp\left[
\beta n\ln n\right]  $ or less$.$ Then our choice of $\beta$ and $\lambda$
\ in (\ref{d0}) gives the expected number $N^{\prime}=o(N)$ of non-saturated
caps. Thus, non-saturated caps typically form a vanishing fraction of not only
$\mu$-caps but also of $\rho$-caps.

5. Next, we proceed with saturated $\mu$-caps and count all those centers
$\mathbf{u}^{\prime\prime}\in$ \textsf{Cov}$(\varepsilon)$ that are left
uncovered by random $\rho$-caps$.$ All caps $C(\mu,\mathbf{z}^{\prime\prime})$
that contain uncovered centers $\mathbf{u}^{\prime\prime}$ are called porous.
For a given $s,$ we show in Lemma \ref{lm:good} that the set $\left\{
\mathbf{u}^{\prime\prime}\right\}  $ forms a very small portion of
\textsf{Cov}$(\varepsilon)$ that has the expected order $\omega^{s+1}<
\exp\left[  -n\ln^{2}n/(3\ln\ln n)\right]  .$ Note that the quantity
\[
\left\vert \text{\textsf{Cov}}(\varepsilon)\right\vert \leq{\vartheta}\left(
S_{r}^{n}\right)  /{\Omega_{\varepsilon}}%
\]
exceeds $N$ by the factor $\Omega_{d}/\Omega_{\varepsilon}$ that only grows as
$\exp\left[  n\ln n\right]  $. Therefore, $\ $the expected size of $\left\{
\mathbf{u}^{\prime\prime}\right\}  $ is $N^{\prime\prime}=o(N).$

6. Finally, the centers of \ all \ non-saturated and porous caps are combined
into the set $\mathbf{\bar{z}=}\{\mathbf{z}^{\prime},$ $\mathbf{z}%
^{\prime\prime}\}.$ Then the set $\{\mathbf{x}\}=\{\mathbf{y,\bar{z}}\}$
completely covers the set \textsf{Cov}$(\varepsilon)$ with the caps
$C(\rho,\mathbf{x}).$ Therefore, $\{\mathbf{x}\}$ also covers $S_{r}^{n}$ with
unit caps. \medskip

\textit{Main proofs. }To prove Theorem \ref{th:new}, we first observe (by
numerical comparison) that the existing bound (\ref{est0}) is tighter for
$n\leq100$ than bound (\ref{d-sph}) of Theorem \ref{th:new}. For this reason,
we shall only consider dimensions $n\geq100$. The proof is based on three lemmas.

Consider two caps $C(\mu,\mathbf{Z})$ and $C(\rho,\mathbf{Y})$ with centers
$\mathbf{Y}$ and $\mathbf{Z,}$ which are $d$\textit{-}close. These caps are
represented in Fig. 1. Here the origin $\mathbf{O}$ is the center of
$S_{r}^{n}$. \begin{figure}[tbh]
\begin{center}
\includegraphics[width=3.5in]{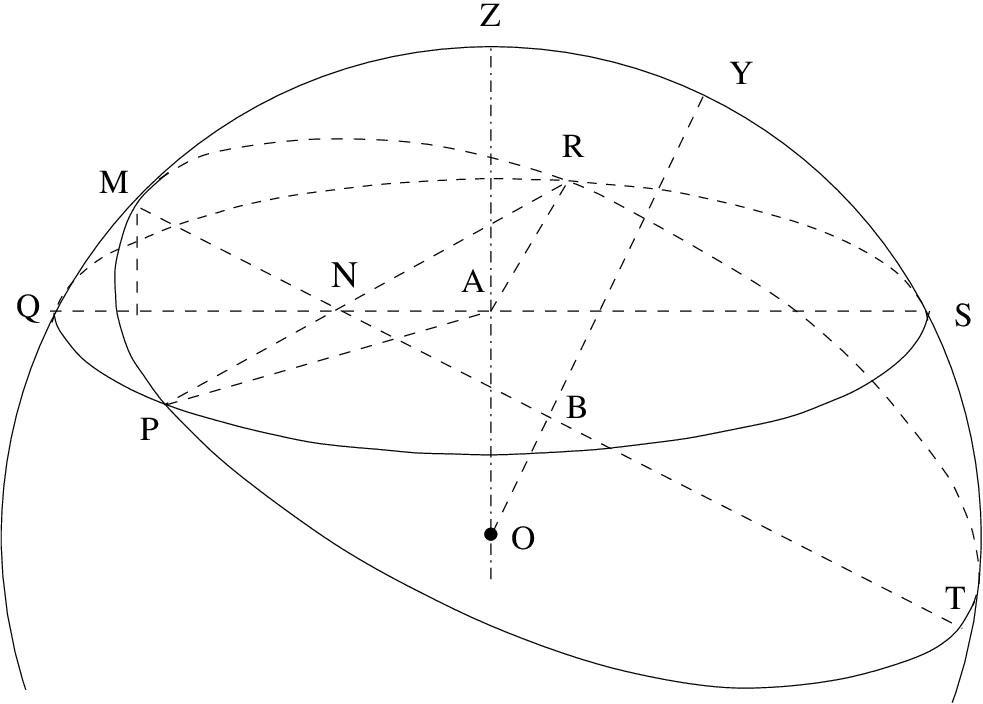}
\end{center}
\caption{Two intersecting caps $C(\mu,\mathbf{Z})$ and $C(\rho,\mathbf{Y})$
with bases $PQRS$ and $PMRT.$}%
\end{figure}

\begin{lemma}
\label{lm:lower}For any cap $C(\mu,\mathbf{Z}),$ a randomly chosen $d$-close
cap $C(\rho,\mathbf{Y})$ fails to cover any given point $\mathbf{x}$ of
$C(\mu,\mathbf{Z})$ with probability $p(\mathbf{x})\leq\omega,$ where
\begin{align}
\omega &  \leq\frac{1}{4\ln n}\exp\{\left(  n-1\right)  \nu_{n}/2\},\label{v}%
\\
\nu_{n}  &  =\ln\left(  1-\frac{4\ln^{2}n}{n}\right)  <-\frac{4\ln^{2}n}%
{n}\nonumber
\end{align}

\end{lemma}

\noindent\textit{Proof.} \ The caps $C(\mu,\mathbf{Z})$ and $C(\rho
,\mathbf{Y})$ have bases\textbf{\ }$\mathbf{PQRSA}$ and $\mathbf{PMRTB},$
which form the balls\ $B_{\mu}^{n}(\mathbf{A})$ and $B_{\rho}^{n}%
(\mathbf{B}).$ Below we consider the boundary $\mathbf{PQRS}$ of
$C(\mu,\mathbf{Z}),$ which forms the sphere $S_{\mu}^{n-1}(\mathbf{A}).$ The
bigger cap $C(\rho,\mathbf{Y})$ covers this boundary, with the exception of
the cap $\mathbf{PQR}$ centered at $Q.$ We first consider the case, when
$\mathbf{x}$ is a boundary point and belongs to $\mathbf{PQRS}.$ Then the
probability $p(\mathbf{x})$ is the fraction $\Omega$ of the entire boundary
contained in uncovered cap $\mathbf{PQR}.$ We first estimate the half-angle
$\alpha=\measuredangle\mathbf{PAQ}$ formed by the cap $\mathbf{PQR}.$

Let $d(\mathbf{H,G})$ denote the distance between any two points $\mathbf{H}$
and $\mathbf{G.}$ Also, let $\sigma(\mathbf{H})$ be the distance from a point
$\mathbf{H}$ to the line $\mathbf{OBY}$ that connects the origin $\mathbf{O}$
with the center $\mathbf{B}$ of the bigger base $B_{\rho}^{n}(\mathbf{B})$ and
with the center $\mathbf{Y}$ of the cap $C(\rho,\mathbf{Y}).$ We use
inequalities
\[
\sigma(\mathbf{A})\leq\sigma(\mathbf{Z})\leq d(\mathbf{Z,Y)\leq}d.
\]
On the other hand, consider the base $\mathbf{PNR}$ of the uncovered cap
$\mathbf{PQR.}$ Here $\mathbf{N}$ denotes the center of this base. Then both
lines $\mathbf{AN}$ and $\mathbf{BN}$ are orthogonal to this base. Also,
$d(\mathbf{B,P})=\rho,$ and $d(\mathbf{N,P})\leq d(\mathbf{A,P})=\mu.$ Thus,
\begin{equation}
d(\mathbf{B,N})=\sqrt{d^{2}(\mathbf{B,P})-d^{2}(\mathbf{N,P})}\geq\sqrt
{\rho^{2}-\mu^{2}}\geq\rho-\mu^{2}/(2\rho)\geq\rho-\mu^{2}.\label{d2}%
\end{equation}
Finally, the center line $\mathbf{OBY}$ is orthogonal to the entire base
$\mathbf{PMRTB}$ and its line $\mathbf{BN.}$ Then $d(\mathbf{B,N}%
)=\sigma(\mathbf{N})$ and%
\begin{equation}
d(\mathbf{A,N})\geq\sigma(\mathbf{N})-\sigma(\mathbf{A})\geq\rho-\mu
^{2}-d=\varepsilon.\label{d3}%
\end{equation}
Now we consider the right triangle $\mathbf{ANP}$ and deduce that
\begin{align}
\cos\alpha &  =d(\mathbf{A,N})/d(\mathbf{A,P})\geq\varepsilon/\mu=\left(  2\ln
n\right)  /\sqrt{n}\label{cos}\\
\sin^{2}\alpha &  \leq\ln\left(  1-\frac{4\ln^{2}n}{n}\right)  =\exp\{-\nu
_{n}\}.\nonumber
\end{align}
We can now estimate the fraction $\Omega$ of the boundary sphere $S_{\mu
}^{n-1}(\mathbf{A})$ contained in the uncovered cap $\mathbf{PQR}.$ Here we
use the bound (see \cite{5}, Corollary 3.2) immediately follows from
\[
\Omega<\left(  \sin^{n-1}\alpha\right)  /\left(  \sqrt{2\pi(n-1)}\cos
\alpha\right)
\]
To obtain (\ref{v}), we simply replace denominator with a smaller quantity
$4\ln n$ using (\ref{cos}).

Finally, consider the second case, when a point $\mathbf{x}$ does not belong
to the boundary $\mathbf{PQRS}.$ Therefore, $\mathbf{x}$ is \ taken from a
smaller cap $C(\mu^{\prime},\mathbf{Z})\subset C(\mu,\mathbf{Z})$ with the
same center $\mathbf{Z}$ and a half-chord $\mu^{\prime}<\mu.$ Similarly to the
first case, we define the boundary $S_{\mu^{\prime}}^{n-1}(\mathbf{A}^{\prime
})$ of \ $C(\mu^{\prime},\mathbf{Z}),$ where $\mathbf{A}^{\prime}$ is some
center on the line $\mathbf{AZ}$. \ Again, $p(\mathbf{x})$ is the fraction
$\Omega^{\prime}$ of this boundary left uncovered by the bigger cap
$C(\rho,\mathbf{Y}).$ To obtain the upper bound on $\Omega^{\prime},$ we only
need to replace $\mu$ with a smaller $\mu^{\prime}$ in (\ref{cos}). \ This
gives the smaller angle $\alpha^{\prime}\leq\arccos(\varepsilon/\mu^{\prime
}),$ which is reduced to $0$ if $\varepsilon\geq\mu^{\prime}.$ In particular,
the center $\mathbf{Z}$ of the cap is always covered by any $d$-close cap$.$
Thus, we see that any internal layer $S_{\mu^{\prime}}^{n-1}(\mathbf{A}%
^{\prime})$ of the cap $C(\mu,\mathbf{Z})$ has a smaller uncovered fraction
$\Omega^{\prime}\leq\Omega.$ This gives the required condition $p(\mathbf{x}%
)\leq\Omega^{\prime}\leq\Omega<\omega$ for any point $\mathbf{x}$ and proves
our lemma. \ {\hfill}\ ${\square}$\smallskip

\textit{Remarks. } First, note that (\ref{cos}) can be used only for
$n\geq75.$ For $n<75,$ inequality (\ref{cos}) gives $\cos\alpha\geq1,$ which
only shows that $\alpha=0$. In this case, a $d$-close cap entirely covers
$\mu$-cap in Fig. 1. Second, note that even a marginal increase in $d$
completely changes our setting. Namely, it can be proven that about half \ the
base of the $\mu$-cap is uncovered if a $\rho$-cap\ is $\left(  d+\varepsilon
\right)  $-close.\smallskip

Our next goal is to estimate the expected number $N^{\prime}$ \ of
non-saturated caps $C(\mu,\mathbf{z}^{\prime})$ left after $N$ trials$,$ where
$N=\left\lfloor \left(  \lambda n\ln n\right)  /\Omega_{d}\right\rfloor $
according to (\ref{lam1}).

\begin{lemma}
\label{lm:bad}For $n\geq100,$ the number of \ non-saturated caps
$C(\mu,\mathbf{z}^{\prime})$ left after $N$ trials has expectation $N^{\prime
}<e^{-n/4-1}N.$
\end{lemma}

\noindent\textit{Proof.} \ \ Given any center $\mathbf{z},$ a randomly chosen
center $\mathbf{y}$ is $d$-close to $\mathbf{z}$ with the probability
$\Omega_{d}.$ Then the probability to obtain at most $s$ such caps
$C(d,\mathbf{y})$ is
\begin{equation}
P=\sum_{i=0}^{s}\left(  _{\,i}^{N}\right)  \Omega_{d}^{i}(1-\Omega_{d}%
)^{N-i}.\label{prob}%
\end{equation}
First, we use (\ref{expan2}) for $\Omega_{d}$ and take $d=1-2\varepsilon
-\mu^{2}$ in (\ref{d0}). Then
\[
\Omega_{d}\geq\theta_{n}=\frac{1}{2}-\sqrt{\frac{1}{\pi\ln^{2}n}+\frac
{1}{16\pi\ln^{6}n}}%
\]
It is easy to verify that $\theta_{n}\geq\theta_{100}>\theta,$ where
$\theta\equiv0.377.$ For brevity, let $\Omega_{d}=x.$ Then $(1-x)^{N-i}%
\leq(1-x)^{N-s}$. Note also that $Nx=L\geq\lambda n\ln n-1/2$ according to
(\ref{lam1}). Then we use Taylor series $\ln\left(  1-x\right)  =-x-x^{2}%
/2-...$ and see that
\begin{align*}
\ln(1-x)^{N-s} &  =-L-x(L/2-s)-x^{2}(L/3-s/2)-...\\
&  =L\left[  \ln(1-x)\right]  /x-s\ln(1-x)
\end{align*}
Every term $x^{i}$ in the first line has negative coefficient $L/(i+1)-s/i$
since $s=n/q.$ Thus, we can take $x=\theta,$ which gives the upper bound
\begin{align}
\ln(1-x)^{N-s} &  \leq L\left[  \ln(1-\theta)\right]  /\theta-s\ln
(1-\theta)\label{prob-a}\\
&  \leq-n\left[  \tau_{1}\lambda\ln n-\tau_{2}/\left(  6\ln\ln n\right)
\right]  +1\nonumber
\end{align}
with coefficients
\[
\tau_{1}=\left[  \ln(1-\theta)\right]  /\theta=1.255,\tau_{2}=-\ln
(1-\theta)=0.473
\]
Now we proceed with the remaining terms in (\ref{prob}).\ Here $\left(
_{\,i}^{N}\right)  x^{i}<\left(  Nx\right)  ^{i}/i!$ $\leq(\lambda n\ln
n)^{i}i!$ and
\begin{gather*}
\sum\nolimits_{i=0}^{s}\left(  _{\,i}^{N}\right)  x^{i}\leq\frac{(\lambda n\ln
n)^{s}}{s!}\sum\nolimits_{i=0}^{\infty}(\lambda q\ln n)^{-i}\leq
3\frac{(\lambda n\ln n)^{s}}{s!}\quad\quad\\
\leq\frac{(\lambda n\ln n)^{s}}{(s/e)^{s}}=(e\lambda q\ln n)^{n/q}%
=\exp\left\{  nh_{n}\right\}  ,\;h_{n}=\ln(e\lambda q\ln n)/q
\end{gather*}
Direct calculations give the bound $\sum\nolimits_{i=0}^{\infty}(\lambda q\ln
n)^{-i}<3$ for any $n\geq100.$ Then we use Sterling formula $s!$ $>(2\pi
s)^{1/2}(s/e)^{s}.$ Summarizing, we have%
\[
P\leq\exp\left\{  n[h_{n}-\tau_{1}\lambda\ln n+\tau_{2}/\left(  6\ln\ln
n\right)  \right\}  ,
\]
Next, we calculate the size $\left\vert \text{\textsf{Cov}}(\mu)\right\vert $
needed to cover $S_{r}^{n}.$ Comparing parameters $\vartheta_{\ast}$ of
(\ref{d-sphe}) and $\lambda$ of (\ref{d0}), we see that $\vartheta_{\ast}%
\leq2\lambda n\ln n-1\leq2\Omega_{d}N$ for any $n\geq100.$Thus, we can cover
$S_{r}^{n}$ with density $\vartheta_{\ast}$ which gives
\begin{equation}
\left\vert \text{\textsf{Cov}}(\mu)\right\vert \leq\vartheta_{\ast}%
/\Omega_{\mu}\leq2N\Omega_{d}/\Omega_{\mu}\label{d}%
\end{equation}
Finally, we use (\ref{p-a}), which gives $\Omega_{\mu}\geq\Omega_{1}\mu
^{n}/\left(  8\sqrt{n}\right)  $ for $n\geq100$ and
\[
\left\vert \text{\textsf{Cov}}(\mu)\right\vert \leq16\sqrt{n}N\mu^{-n}%
\leq16\sqrt{n}N\exp\left\{  \beta n\ln n+n\ln4\right\}
\]
Thus, the expected number of \ non-saturated caps is
\begin{align}
N^{\prime} &  \leq\left\vert \text{\textsf{Cov}}(\mu)\right\vert P\leq
16\sqrt{n}N\exp\left\{  n\Psi_{n}\right\}  \nonumber\\
\Psi_{n} &  =\ln(e\lambda q\ln n)/q-\left(  \tau_{1}\lambda-\beta\right)  \ln
n+\tau_{2}/q+\ln4+1/n\smallskip\label{bad1}%
\end{align}
The quantity $\Psi_{n}$ has moderating term $-\left(  \tau_{1}\lambda
-\beta\right)  \ln n$ and declines in $n.$ Direct calculation shows that
$\Psi_{_{100}}<-0.36.$ Finally, we replace $16\sqrt{n}\exp\left\{  -n\Psi
_{n}\right\}  $ with a larger quantity $\exp\left\{  -n/4-1\right\}  .$%
{\hfill\ }${\square}$\smallskip

Consider now the saturated caps $C(\mu,\mathbf{z})$ and the centers
$\mathbf{u\in}$\textsf{Cov}$(\varepsilon)$ inside them.

\begin{lemma}
\label{lm:good}For any $n\geq100,$ the number of centers $\mathbf{u}%
^{\prime\prime}\in$\textsf{Cov}$(\varepsilon)$ left uncovered in all saturated
caps $C(\mu,\mathbf{z})$ has expectation $N^{\prime\prime}<e^{-2n}N.$
\end{lemma}

\noindent\textit{Proof}. We first estimate the total number $\left\vert
\text{\textsf{Cov}}(\varepsilon)\right\vert $ of centers $\mathbf{u}.$ We
proceed similarly to (\ref{d}). Then
\[
\left\vert \text{\textsf{Cov}}(\varepsilon)\right\vert \leq2N\Omega_{d}%
/\Omega_{\varepsilon}\leq16\sqrt{n}N\left(  2n\ln^{2}n\right)  ^{n}.
\]
Any cap $C(\mu,\mathbf{z})$ intersects with at least $\ s+1$ randomly chosen
caps $C(\rho,\mathbf{y}).$ According to Lemma {\ }\ref{lm:lower}, any single
$\rho$-cap fails to cover any given point $x\in C(\mu,\mathbf{z})$ with
probability $\omega$ or less. Therefore any point $\mathbf{u}^{\prime\prime
}\in C(\mu,\mathbf{z})$ is not covered with probability at most $\omega
^{s+1}<\omega^{n/q}\leq\exp\left\{  nC_{n}\right\}  ,$ where we use (\ref{v})
and obtain
\[
C_{n}=\frac{n-1}{12\ln\ln n}\ln\left(  1-\frac{4\ln^{2}n}{n}\right)
-\frac{\ln(4\ln n)}{6\ln\ln n}%
\]
$\ $ Then
\begin{align}
N^{\prime\prime}  &  \leq\left\vert \text{\textsf{Cov}}(\varepsilon
)\right\vert \cdot\omega^{n/q}\leq16N\sqrt{n}\exp\left\{  n\Phi_{n}\right\}
,\label{p}\\
\Phi_{n}  &  =C_{n}+\ln n+2\ln\ln n+\ln2.\nonumber
\end{align}
Note that the first term in $C_{n}$ has the order of $-\left(  \ln
^{2}n\right)  /\left(  3\ln\ln n\right)  $ and moderates all other terms.
Thus, $\Phi_{n}<\Phi_{_{100}}.$ Direct calculation shows that $\Phi_{_{100}%
}<-2,$ which proves the lemma.\ \ {\hfill}\ ${\square}$\medskip

\noindent\textit{Proof of Theorem \ref{th:new}.} Consider any cap
$C(\mu,\mathbf{\bar{z}})$ that contains at least one uncovered center
$\mathbf{u}\in\,$\textsf{Cov}$(\varepsilon).$ Such a cap is either
non-saturated or porous and therefore $\{\mathbf{\bar{z}}\}=\left\{
\mathbf{z}^{\prime}\right\}  \cup\left\{  \mathbf{z}^{\prime\prime}\right\}
.$ Then, according to Lemmas \ref{lm:bad} and \ref{lm:good}, $\{\mathbf{\bar
{z}}\}$ has expected size $\bar{N}\leq N^{\prime}+N^{\prime\prime}<e^{-n/4}N.$
Thus, there exist $N$ randomly chosen centers $\mathbf{y}$ that leave at most
$e^{-n/4}N$ centers $\mathbf{\bar{z}.}$ The extended set $\ \left\{
\mathbf{x}\right\}  =\left\{  \mathbf{y}\right\}  \cup\{\mathbf{\bar{z}}\}$
forms a unit covering of $S_{r}^{n}.$ This covering has density%
\[
\vartheta\leq\Omega_{1}\left(  N+\bar{N}\right)  \leq\Omega_{1}N(1+e^{-n/4}%
)\leq\lambda n\ln n\left(  1+e^{-n/4}\right)  \Omega_{1}/\Omega_{d}%
\]
Similarly to inequality (\ref{expan2}), we can directly verify that
$\Omega_{1}/\Omega_{d}\leq1+3/\left(  2\ln n\right)  $ for $n\geq100.$
Finally, we take $\lambda$ of (\ref{d0}) and combine the last inequalities for
$\vartheta$ and $\Omega_{d}/\Omega_{1}$ as follows
\begin{align*}
\frac{\vartheta}{n\ln n} &  \leq\left(  \frac{1}{2}+\frac{3\ln\ln n}{\ln
n}+\frac{1}{2\ln n}\right)  \left(  1+\frac{3}{2\ln n}\right)  (1+e^{-n/4})\\
&  <\frac{1}{2}+\frac{3\ln\ln n}{\ln n}+\frac{3}{\ln n}.
\end{align*}
Direct verification shows that the last inequality holds for $n=100.$ Then we
can again use Lemma \ref{lm:mod} for larger $n.$ This completes the proof of
Theorem \ref{th:new}. \ {\hfill}\ ${\square}$\smallskip

Finally, note that Theorem \ref{th:new} directly leads to Corollary
\ref{th:ball}. Indeed, we can use the well known fact $\vartheta
(\mathbb{R}^{n-1})=\lim_{r\rightarrow\infty}\vartheta(S_{r}^{n})$ (see the
proof in \cite{2} or Theorem $\Pi.1$ in \cite{lev}, where a similar proof is
detailed for packings of $\mathbb{R}^{n}$). \smallskip

\emph{Further directions.} We have proved that the classical Rogers bound
(\ref{d-ball}) on the covering density of a sphere $S_{r}^{n}$ or the
Euclidean space $\mathbb{R}^{n}$ can be reduced about two times for large
dimensions $n.$ One open problem is to reduce the gap between this bound and
its lower counterpart (\ref{few}), which is linear in $n.$ In this regard,
note that our design always holds if parameter $\mu$ has the order of
$n^{-\beta}$ given some constant $\beta>1/2$. However, it can be verified that
choosing a smaller constant $\beta<1/2$ will exponentially increase the
covering size. Therefore, a completely new design is needed for any further
reductions of the upper bound. Another important problem is to extend the
above results to the balls $B_{r}^{n}$ of an arbitrary radius $r.$ Our
conjecture is that $\vartheta(B_{r}^{n})\leq\left(  1/2+o(1)\right)  n\ln n$
for any $r$ and $n\rightarrow\infty.$\smallskip

\textit{Addendum. }The published version of this paper \cite{dum2} uses
incorrect formula (8) $\Omega_{\tau}/\Omega_{\rho}\geq\left(  \tau
/\rho\right)  ^{n}$ instead of the correct inequality $\Omega_{\tau}%
/\Omega_{\rho}\geq\left(  \arcsin\tau\,/\,\arcsin\rho\right)  ^{n}$ valid for
any $\tau<\rho\leq1$. Formula (8) was not used in main Lemmas 6,7 and 8 of
paper \cite{dum2}; however, it reduced the residual terms in the preliminary
Section 2. The current version excludes formula (8) and also tightens several
estimates of \cite{dum2}. The author wishes to thank Dr. Marton Naszodi, who
notified the author of this error.

\end{document}